\documentclass[12pt, reqno]{amsart}
\usepackage{amsmath, amsthm, amscd, amsfonts, amssymb, graphicx, color}
\usepackage[bookmarksnumbered, colorlinks, plainpages]{hyperref}
\input{mathrsfs.sty}

\hypersetup{colorlinks=true,linkcolor=red, anchorcolor=green,
citecolor=cyan, urlcolor=red, filecolor=magenta, pdftoolbar=true}

\newtheorem{twr}{Theorem}[section]
\newtheorem{lem}[twr]{Lemma}

\newtheorem{cor}[twr]{Corollary}
\theoremstyle{definition}

\theoremstyle{remark}

\numberwithin{equation}{section}

\DeclareMathOperator{\conv}{conv}

\DeclareMathOperator{\lin}{lin}

\DeclareMathOperator{\tr}{Tr}

\DeclareMathOperator{\id}{Id}
\DeclareMathOperator{\dist}{dist}

\DeclareMathOperator{\ext}{ext}
\DeclareMathOperator{\Gr}{Gr}
\usepackage{color}
\usepackage{enumerate}

\begin{document}

\setcounter{page}{1}

\title[On the dimension of the set of minimal projections]
{On the dimension of the set of minimal projections}
\author[T. Kobos \MakeLowercase{and} G. Lewicki]
{Tomasz Kobos$^1$ \MakeLowercase{and} Grzegorz Lewicki$^2$}
\address{$^1$Faculty of Mathematics and Computer Science, Jagiellonian University, \L ojasiewicza 6, 30-348 Krak\'ow, Poland}
\email{tomasz.kobos@uj.edu.pl}
\address{$^2$Faculty of Mathematics and Computer Science, Jagiellonian University, \L ojasiewicza 6, 30-348 Krak\'ow, Poland}
\email{grzegorz.lewicki@uj.edu.pl}

%%%%%%%%%%%%%%%%%%%%%%%%%%%%%%%%%%%%%%%%%%%%%%
\begin{abstract}
Let $X$ be a finite-dimensional normed space and let $Y \subseteq X$ be its proper linear subspace. The set of all minimal projections from $X$ to $Y$ is a convex subset of the space all linear operators from $X$ to $X$ and we can consider its affine dimension. We establish several results on the possible values of this dimension. We prove optimal upper bounds in terms of the dimensions of $X$ and $Y$. Moreover, we improve these estimates in the polyhedral normed spaces for an open and dense subset of subspaces of the given dimension. As a consequence, in the polyhedral normed spaces a minimal projection is unique for an open and dense subset of hyperplanes. To prove this, we establish certain new properties of the Chalmers-Metcalf operator. Another consequence is the fact, that for every subspace of a polyhedral normed space, there exists a minimal projection with many norming pairs.
\end{abstract} \maketitle
%%%%%%%%%%%%%%%%%%%%%%%%%%%
\section{Introduction}

If $X$ is a normed space over $\mathbb{K}$ (where $\mathbb{K}=\mathbb{R}$ or $\mathbb{K}=\mathbb{C}$) and $Y \subseteq X$ is its proper linear subspace, then by \emph{projection} from $X$ to $Y$ we shall mean a bounded linear operator $P:X \to Y$ such that $P|_{Y}=\id_{Y}$. By $\mathcal{P}(X, Y)$ we denote the set of all projections from $X$ onto $Y$. The \emph{relative projection constant} of $Y$ is defined as
$$\lambda(Y, X) = \inf \{ \|P\| : P \in \mathcal{P}(X, Y) \}.$$
Moreover, if a projection $P: X \to Y$ satisfies $\|P\| = \lambda(Y, X)$ then $P$ is called a \emph{minimal projection}. The set of all minimal projections will be denoted by $\mathcal{P}_{\min}(X, Y)$. By $B_X$ and $S_X$ we will always denote the unit ball and the unit sphere of $X$, respectively.

Projection is one of the fundamental concepts of the functional analysis and minimal projections have been studied extensively for many years. The following questions are classical, when considered for some specific normed spaces $Y \subseteq X$ (see for example: \cite{basso}, \cite{blatter}, \cite{cheneyfourier}, \cite{chalmersmetcalf}, \cite{CM1}, \cite{CM2}, \cite{bb}, \cite{bb2} \cite{grunbaum}, \cite{konig}, \cite{lewickichalm}, \cite{lewickichalm2}, \cite{lewickitwodim} \cite{lewickiskrzypeklp}, \cite{lozinski}, \cite{odinec}, \cite{odinec2}, \cite{skrzypektwodim}, \cite{skrzypeknonunique}, \cite{skrzypekmatrix}, \cite{sokolowski}):
\begin{enumerate}
    \item Determine $\lambda(Y, X)$.
    \item Determine a minimal projection $P_0: X \to Y$, that is some projection $P_0$ satisfying $\|P_0\| = \lambda(Y, X)$.
    \item Determine if such projection $P_0$ is unique.
\end{enumerate}

One famous example is Lozinski Theorem about the minimality of the classical Fourier projection (see \cite{lozinski}). The minimality of Fourier projection was proved by Lozinski in 1948, but it took another 20 years to establish, that it is in fact, the unique minimal projection. It was proved by Cheney, Hobby, Morris, Schurer and Wulbert (see \cite{cheneyfourier}). This shows that, generally speaking, even if the second question is settled, the third one can still provide an additional and significant challenge. The main goal of this paper is to consider the third question, but in a much broader sense than traditionally taken. To put our investigation in a context, we will recall some results related to the uniqueness of minimal projection.

An old theorem due to Odyniec (\cite{odinec2}) can be considered as a starting point of our investigation. It states, that if $X$ is a three-dimensional real normed space and $Y \subseteq X$ is a two-dimensional subspace such that $\lambda(Y, X)>1$, then the minimal projection $P:X \to Y$ is unique. Moreover, Odyniec proved a similar result for an $n$-dimensional real normed space $X$ and an $(n-1)$-dimensional subspace $Y \subseteq X$ satisfying $\lambda(Y, X)>1$, but with the additional assumption of the uniform convexity (\cite{odinec}). On the other hand, for the $1$-complemented subspaces we have a classical result of Cohen and Sulivan (\cite{cohen}), who proved that for $1$-complemented subspace $Y$ in a smooth Banach space $X$, a minimal projection is unique. This is complemented by a result of Shekhtman and Skrzypek (\cite{skrzypektwodim}), who proved that if $X$ is both uniformly convex and smooth, then there exists a unique minimal projection onto every two-dimensional subspace of $X$. However, later in \cite{skrzypeknonunique} the same authors provided an example of a three-dimensional subspace of $\ell_p^5$ (where $1<p<\infty$), for which a minimal projection is not unique. Thus, smoothness and uniform convexity are not enough to guarantee the uniqueness of the minimal projection, if the subspace is of dimension greater than $2$ and codimension greater than $1$.

If a minimal projection turns out to be not unique, then the third question is answered, but it is still possible to pursuit it further.  One can consider the set of all minimal projections $\mathcal{P}_{\min}(X, Y)$ and try to grasp how "large" is it. This could be defined in a number of ways, but one of the possibilities, is to consider its dimension. It is a straightforward observation that the set $\mathcal{P}_{\min}(X, Y)$ is a convex subset in the space of all linear operators $\mathcal{L}(X, X)$. Thus, one can consider its \emph{affine dimension}, defined as the minimal possible dimension of an affine subspace $V \subseteq X$ such that $\mathcal{P}_{\min}(X, Y) \subseteq V$. When dealing with convex subsets of real normed spaces, we will refer to this as simply dimension and we will denote it as $\dim$. Our aim is to gain some insight about $\dim \mathcal{P}_{\min}(X, Y)$.

If a minimal projection is unique, then $\dim \mathcal{P}_{\min}(X, Y)=0$, but in general, this number can achieve many different values. In this paper paper, we study the possible values of $\dim \mathcal{P}_{\min}(X, Y)$, when $X$ is a finite-dimensional normed space (mostly real) and $Y \subseteq X$ is a proper linear subspace. Interestingly, just as in the aforementioned theorem of Odyniec, in our results there often seems to be an important distinction between the case of $1$-complemented subspaces and not $1$-complemented subspaces.

Let us suppose that $X$ is an real $n$-dimensional normed space and $Y \subseteq X$ is a linear subspace of dimension $1 \leq k \leq n-1$. In Section \ref{sectgeneral} we establish the following general upper bounds for $\dim \mathcal{P}_{\min}(X, Y)$ (see Theorem \ref{twrgeneral}):
\begin{enumerate}
    \item $\dim \mathcal{P}_{\min}(X, Y) \leq k(n-k)$.
    \item If $\lambda(Y, X)>1$, then $\dim \mathcal{P}_{\min}(X, Y) \leq k(n-k)-2$.
\end{enumerate}
Both estimates are optimal. Thus, the second estimate is a broad generalization of the theorem of Odyniec, which corresponds to the case of $n=3$ and $k=2$.

Even if these estimates are optimal, one can ask to what degree these estimates reflect what happens on a "regular basis". This is a broad question, to which we can provide a partial answer in a specific class of norms. We say that the normed space $X=(\mathbb{R}^n, \| \cdot \|)$ is \emph{polyhedral}, if the unit ball of $X$ is a convex polytope in $\mathbb{R}^n$, i.e. $B_X = \conv \{ \pm x_1, \pm x_2, \ldots, \pm x_N \}$ for some $x_1, x_2, \ldots, x_N \in \mathbb{R}^n$. The importance of polyhedral normed spaces is a direct consequence of the importance of the class of convex polytopes when considered as a subset of all convex bodies. It turns out, that our estimates can be significantly improved for polyhedral normed spaces and almost all subspaces $Y$. Let us suppose that $X$ is an $n$-dimensional polyhedral normed space and $Y \subseteq X$ is a linear subspace of dimension $1 \leq k \leq n-1$. If $Y$ is in a general position (for the precise meaning we refer to Section \ref{sectpoly}), then
\begin{equation}
\label{estimatepoly}
\dim \mathcal{P}_{\min}(X, Y) \leq k(n-k)-n+1.
\end{equation}
This estimate holds for an open and dense subset of $k$-dimensional subspaces $Y$. In particular, for $k=n-1$, we obtain that the hyperplane projections are unique in the polyhedral norms for an open and dense set of hyperplanes. This generalizes some previously known results, which have dealt with some particular norms. In these results, one could observe, that a minimal projection on the hyperplane is not unique only in some quite specific situations. See Theorem \ref{twrwymiar} and whole Section \ref{sectpoly} for a broader discussion. We do not know, if for $2 \leq k \leq n-2$ this estimate is optimal.

Throughout the paper we rely a lot on the Chalmers-Metcalf operator, which has already proved to be a highly effective tool, when dealing with the minimal projections. It was first used by Chalmers and Metcalf in \cite{chalmersmetcalf} to find a minimal projection from $C[-1, 1]$ onto the linear space of polynomials of degree at most $2$. It is worth mentioning, that to prove estimate (\ref{estimatepoly}) we establish a property of Chalmers-Metcalf operator, which is a broad generalization of the one that was commonly used before (see Lemma \ref{lemcmdim}). As a by-product of our investigation, we establish, that in the polyhedral normed spaces, for every subspace there exists a minimal projection with at least $n$ norming pairs (see Theorem \ref{twrpunkty}).

All required properties of the Chalmers-Metcalf operators and other preliminaries from the general approximation theory are introduced in Section \ref{sectapprox}. It should be strongly emphasized, that throughout the paper, we will rely heavily on the background information and notation presented in this section. We will also refer often to the $\ell_1^n$ and $\ell_{\infty}^n$ spaces, which will be identified with $\mathbb{K}^n$ with the respective norm (denoted by $\| \cdot \|_1$ and $\| \cdot \|_{\infty}$). For a subset $A$ of a normed space, by $\lin A$ we will denote the linear span of $A$.

\section{Preliminaries from the general approximation theory and the Chalmers-Metcalf operator}
\label{sectapprox}
%definicja supported
%zbior najlepszych aproksymacji jest wypukly

We start the discussion with a more general setting than considered in the rest of the paper. The classical problem of the approximation theory is to describe the set 
$$P_Y(x)=\{y \in Y: \ \|x-y\|=\dist(x, Y)\}$$
of best approximations, where $x \in X$ is an element of a normed space $X$ (or, even more generally, a metric space) and $Y \subseteq X$ is a non-empty subset. To simplify the discussion, we will stick to the setting of $X$ being a finite-dimensional normed space and $Y \subseteq X$ being a proper linear subspace. Many of the discussed properties could be considered also in the infinite-dimensional setting, but as the rest of the paper is concerned only with the finite-dimensional case, we will omit general considerations. The problem of determining the set of minimal projections can be related to the question of best approximation in the following way. If $X$ is a $n$-dimensional space over $\mathbb{K}$ and $Y \subseteq X$ is a linear subspace of dimension $1 \leq k \leq n-1$, then we can consider a normed space $\mathcal{L}_Y(X, Y) \subseteq \mathcal{L}(X, Y)$, which consists of all linear operators $T:X \to Y$ satisfying $T|_Y \equiv 0$. Clearly, this is a normed space with the usual operator norm and if we pick any basis $y_1, y_2, \ldots, y_k$ of $Y$ and any linearly independent set of functionals $f_1, f_2, \ldots f_{n-k} \in Y^{\perp}$ (where $Y^{\perp}$ denotes a linear subspace of $X^*$ consisting of functionals vanishing on $Y$), then the set
$$\{y_i \otimes f_j: \ 1 \leq i \leq k, 1 \leq j \leq n-k \} \subseteq \mathcal{L}_Y(X, Y),$$
forms a basis of the space $\mathcal{L}_Y(X, Y)$. Here, if $x_0 \in X$ and $f_0 \in X^*$, then by $x \otimes f \in \mathcal{L}(X, X)$ we understand a linear operator of rank one: 
$$X \ni x \to f_0(x)x_0 \in X.$$
In particular, we have $\dim \mathcal{L}_Y(X, Y) = k(n-k)$. It is now evident that:
\begin{itemize}
    \item $P_0-P_1 \in \mathcal{L}_Y(X, Y)$ for any two projections $P_0, P_1 \in \mathcal{P}(X, Y)$,
    \item $P_0 \in \mathcal{P}_{\min}(X, Y)$ if and only if, $0 \in P_{\mathcal{L}_Y(X, Y)}(P_0)$,
    \item For any projection $P_0 \in \mathcal{P}(X, Y)$ we have $\mathcal{P}_{\min}(X, Y) = P_{\mathcal{L}_Y(X, Y)}(P_0) + P_0$.
\end{itemize}
Thus, the question of determining the set of minimal projection is equivalent to the classical question of finding the set of best approximations, considered for any fixed projection $P_0 \in \mathcal{P}(X, Y)$ in the linear subspace $\mathcal{L}_Y(X, Y)$. What will be especially important for us, is a straightforward observation, that the dimension of the set $\mathcal{P}_{\min}(X, Y)$ is the largest $d$, for which one can find linearly independent operators $L_1, L_2, \ldots, L_d \in \mathcal{L}_Y(X, Y)$ and a projection $P_0 \in \mathcal{P}(X, Y)$ such that:
$$P_0, P_0 + L_1, \ldots, P_0 + L_d \in \mathcal{P}_{\min}(X, Y).$$
In the general setting of $X$ being a finite dimensional normed space and $Y \subseteq X$ being a linear subspace, there is a well-known characterization of the set of best approximations using linear functionals.
\begin{twr}
\label{twrfunkcjonal}
Let $X$ be a finite-dimensional normed space and let $Y \subseteq X$ be a proper linear subspace. If $x_0 \in X \setminus Y$, then $y_0 \in P_Y(x)$ if and only if, there exists a functional $f \in B_{X^*}$ such that $f(x-y_0)=\|x-y_0\|$ and $f|_{Y} \equiv 0$. 
\end{twr}

The equality $f(x-y_0)=\|x-y_0\|$ could be also stated as $f(x)=\dist(x, Y)$. Furthermore, any functional $f \in B_{X^*}$ can be written as a convex combination of some extreme points of $B_{X^*}$. Hence, for the $f$ as above, we can write
$$f=\sum_{i=1}^{l} \alpha_if_i,$$
where $\{f_1, f_2, \ldots f_l\} \subseteq \ext B_{X^*}$ and $\alpha_1, \alpha_2, \ldots, \alpha_l$ are positive real numbers, which satisfy $\sum_{i=1}^{l} \alpha_i = 1$. Because $\|f_i\|=1$ for every $1 \leq i \leq l$, it is easy to see that we must have $f_i(x-y_0)=\|x-y_0\|$ for every $1 \leq i \leq l$.

This motivates the following definition, that was introduced in \cite{sudolski}. If $x_0 \in X \setminus Y$, then a set of functionals $\{f_1, f_2, \ldots f_l\} \subseteq \ext B_{X^*}$ will be called an \emph{$I$-set for $Y$ with respect to $x_0$}, if there exists $y_0 \in P_Y(x_0)$ and positive real numbers $\alpha_1, \alpha_2, \ldots, \alpha_l$ with the sum $1$, such that:
\begin{equation}
\label{iset1}
f_i(x_0-y_0)=\|x_0-y_0\| \quad \text{ for every } 1 \leq i \leq l,
\end{equation}
\begin{equation}
\label{iset2}
\sum_{i=1}^{l} \alpha_if_i(y) =0 \quad \text{for every } y \in Y.
\end{equation}
It should be noted that if the properties (\ref{iset1}) and (\ref{iset2}) hold for some element $y_0 \in P_Y(x)$, then they hold for every element $y_1 \in P_Y(x)$. In fact, since $y_0 - y_1 \in Y$ we have
$$\sum_{i=1}^{l} \alpha_if_i(x_0 - y_1) = \sum_{i=1}^{l} \alpha_if_i(x_0 - y_0) + \sum_{i=1}^{l} \alpha_if_i(y_0-y_1)$$
$$=\sum_{i=1}^{l} \alpha_if_i(x_0 - y_0) = \|x_0 - y_0\| = \|x_0-y_1\|.$$
It is now clear, that if for some $1 \leq i \leq l$ we would have $f_i(x_0-y_1) < \|x_0-y_1\|$, then by the triangle inequality it would follow also that $\sum_{i=1}^{l} \alpha_if_i(x_0 - y_1) < \|x_0-y_1\|.$ Thus, for every $1 \leq i \leq l$ we have $f_i(x_0-y_1)=\|x_0-y_1\|$.

We further say that an $I$-set is \emph{minimal}, if it does not contain any proper subset, which forms an $I$-set for $Y$ with respect to $x_0$. The following lemma gives an important property of a minimal $I$-set in the real normed space.

\begin{lem}
 Let $X$ be an $n$-dimensional real normed space and let $Y \subseteq X$ be its linear subspace of dimension $k$, where $1 \leq k \leq n-1$. Suppose that $x_0 \in X \setminus Y$ and $\{f_1, f_2, \ldots, f_l\} \subseteq \ext B_{X^*}$ is a minimal $I$-set for $Y$ with respect to $x_0$, where $l \geq 2$. Then, the functionals $f_1|_Y, f_2|_Y, \ldots, f_{l-1}|_Y$ are linearly independent in $Y^*$.
\end{lem}

\emph{Proof}. It is enough to mimick the proof of the Caratheodory's Theorem. We assume that the conditions (\ref{iset1}) and (\ref{iset2}) are satisfied. Let us now assume that there exist $a_1, a_2, \ldots, a_{l-1} \in \mathbb{R}$, not all zero, such that $\sum_{i=1}^{l-1} a_i f_i = 0$. Let $V \subseteq \mathbb{R}^l$ be a linear subspace defined as
$$V = \{(v_1, v_2, \ldots, v_l) \in \mathbb{R}^l \ : \ v_1f_1 + \ldots + v_lf_l = 0 \}.$$

By the assumption we have $(\alpha_1, \ldots, \alpha_{l-1}, \alpha_l), (a_1, \ldots, a_{l-1}, 0) \in V$. These two vectors are clearly linearly independent and hence $\dim V \geq 2$. Thus, there exists $b=(b_1, \ldots, b_l) \in V$, $b \neq 0 $ such that $b_1 + \ldots + b_l = 0$. In particular, we have $b_i > 0$ for some $1 \leq i \leq l$. We denote
$$r = \min_{1 \leq i \leq l} \left \{\frac{\alpha_i}{b_i} \ : \ b_i>0 \right \} = \frac{\alpha_j}{b_j}.$$

If we now define $\alpha'_i = \alpha_i - rb_i$ for $1 \leq i \leq l$, then $\alpha'_i \geq 0$, $\alpha'_j=0$ and $\sum_{i=1}^{l} \alpha'_i = 1.$ Therefore, a functional $F \in X^*$ defined as
 $$F=\alpha'_1 f_1 + \ldots + \alpha'_l f_l$$
 satisfies $F(x_0)=\dist(x_0, Y)$, $F|_Y \equiv 0$ and is supported on at most $l-1$ functionals. This contradicts the minimality of the given $I$-set and the proof is finished. \qed
 
 It turns out, that it is possible to estimate the dimension of the set of best approximations, using the cardinality of a minimal $I$-set ($P_Y(x)$ is a convex set if $Y$ is a linear subspace). The following lemma, generalizes a result of Sudolski and Wójcik (\cite{sudolski}), who proved that if there exists a minimal $I$-set of cardinality $\dim Y + 1$, then $Y$ is a Chebyshev subspace, that is the best approximation is always unique (in fact, they proved a stronger result, that in this case, the best approximation is even $1$-strongly unique). We note that a minimal $I$-set is always of the cardinality at most $\dim Y + 1$, which follows from the previous Lemma.
 
\begin{lem}
\label{lemosz}
 Let $X$ be an $n$-dimensional real normed space and let $Y \subseteq X$ be its $k$-dimensional subspace, where $1 \leq k \leq n-1$. Suppose that $x_0 \in X \setminus Y$ and there exists a minimal $I$-set for $Y$ with respect to $x_0$ of the cardinality $l$. Then 
 $$\dim P_Y(x_0) \leq k-l+1.$$
 \end{lem}
 
 \emph{Proof}. Clearly the dimension of $P_Y(x_0)$ is not greater than of $Y$ itself and the estimate is therefore evident in the case of $l=1$. Thus, we can suppose that $2 \leq l \leq k+1$.
 
 Let us take $y_0 \in P_Y(x_0)$ and denote $x_1=x_0-y_0$. Since $Y$ is a linear subspace, we have that $0 \in P_Y(x_1)$. If we take the linear span $Z = \lin P_Y(x_1) \subseteq Y$, then our goal is to prove that $\dim Z \leq k-l+1$. Suppose that  $\{f_1, f_2, \ldots, f_l\} \subseteq \ext B_{X^*}$ is a minimal $I$-set for $Y$ with respect to $x_0$. By the previous Lemma, the functionals $f_1|_Y, f_2|_Y, \ldots, f_{l-1}|_Y$ are linearly independent in $Y^*$. Hence, a linear subspace subspace $W=\left ( \bigcap_{i=1}^{l-1} \ker f_i \right ) \cap Y$ is of the dimension $k-l+1$. If $y \in P_Y(x_1)$, then from the fact that also $0 \in P_Y(x_1)$ and condition (\ref{iset1}), it follows that
 $$f_i(x_1-y)=\|x_1-y\|=\|x_1\|=f_i(x_1),$$
 for every $1 \leq i \leq l$ (we recall here that properties (\ref{iset1}) and (\ref{iset2}) hold for every element of the best approximation). Hence, we have $f_i(y)=0$ for every $1 \leq i \leq l$ and in particular we have that $y \in W$. This shows that $Z \subseteq W$, which implies $\dim Z \leq \dim W = k-l+1$ and the proof is finished.\qed

Coming back to the setting of linear projections, the functional from Theorem \ref{twrfunkcjonal} yields the \emph{Chalmers-Metcalf functional} for $Y$. By the characterization of the set of extreme points of the unit ball in the space of linear operators and its dual, the $I$-set for Chalmers-Metcalf functional will consist of certain tensor products. To be more specific, if $x \in X$ and $f \in X^*$, then we define $x \otimes f \in \mathcal{L}^*(X, X)$ as $(x \otimes f)(L) = f(L(x))$ for $L: X \to X$. Caution is needed, as it is a different understanding of $x \otimes f$ than one that was introduced before. This double meaning of this operation will be used throughout the paper, but the context will always readily imply the correct interpretation. If $X$ is a finite-dimensional normed space and if $Y \subseteq X$ is a linear subspace, then by the Chalmers-Metcalf functional for $Y$, we will understand the functional of the form
$$T=\alpha_1 x_1 \otimes f_1 + \alpha_2 x_2 \otimes f_2 + \ldots + \alpha_l x_l \otimes f_l \in \mathcal{L}^*(X, X),$$
which satisfies the conditions:
\begin{itemize}
\item $\alpha_i$ are positive reals with sum $1$, 
\item $T|_{\mathcal{L}_Y(X, Y)} \equiv 0$ 
\item $f_i(P_0(x_i))=\|P_0\|=\lambda(Y, X)$ for every $1 \leq i \leq l$ and some fixed minimal projection $P_0 \in \mathcal{P}_{\min}(X, Y)$.
\end{itemize}
A Chalmers-Metcalf functional is exactly the functional from Theorem \ref{twrfunkcjonal}, applied for the space $\mathcal{L}(X, X)$ and its linear subspace $\mathcal{L}_Y(X, Y)$. In this case, the vector $x_0$ could be any minimal projection and the same functional would work, so we just say about a Chalmers-Metcalf functional for $Y$. A Chalmers-Metcalf functional does not have to be unique.

A Chalmers-Metcalf functional has also an operator version, which is more convenient in some instances. This variant uses a notion of the tensor product $x \otimes f$ introduced before, that is, interpreted as a one rank operator. If $X$ is a finite-dimensional normed space and $Y \subseteq X$ is a linear subspace, then by the \emph{Chalmers-Metcalf operator} for $Y$, we will understand the operator $T:X \to X$ of the form
$$T=\alpha_1 x_1 \otimes f_1 + \alpha_2 x_2 \otimes f_2 + \ldots + \alpha_l x_l \otimes f_l,$$
which satisfies the conditions: 
\begin{itemize}
\item $\alpha_i$ are positive reals with sum $1$, 
\item $T(Y) \subseteq Y$
\item $f_i(P_0(x_i))=\|P_0\|=\lambda(Y, X)$ for every $1 \leq i \leq l$ and some fixed minimal projection $P_0 \in \mathcal{P}_{\min}(X, Y)$.
\end{itemize}
There is a correspondence between Chalmers-Metcalf operators and functionals, with exactly the same pairs $(x_i, f_i)$ appearing in the definition, with only the interpretation of $x_i \otimes f_i$ changing. The main difference lies in the fact, that in the functional version we have a condition of vanishing on $\mathcal{L}_Y(X, Y)$, while in the operator version, $Y$ has to be an invariant subspace of $T$. Which variant is more suitable, depends on a particular situation and we will use both of them freely.

We will mention two further properties of the Chalmers-Metcalf functional (or operator) that will be important to us. If $P:X \to Y$ is a linear projection, then $(x, f) \in \ext B_{X} \times \ext B_{X^*}$ will be called a \emph{norming pair} for $P$, if we have $f(P(x))=\|P\|$. It is important to note, that we will consider norming pairs only among the extreme points of the unit ball and its dual. This will be of crucial importance, when we will be working with the polyhedral norms. We have the following two properties:
\begin{itemize}
\item A Chalmers-Metcalf functional/operator for $Y$ does not depend on the choice of a particular minimal projection. In particular, the pairs $(x_i, f_i)$ are norming pairs for \textbf{all} projections in $\mathcal{P}_{\min}(X, Y)$. On the other hand, it is known that minimal projections can possibly have some norming pairs that do not appear in the Chalmers-Metcalf functional/operator.
\item The norm of the minimal projection, that is the relative projection constant of $Y$, can be retrieved easily from the Metcalf-Chalmers operator: we have an equality $\tr (T|_Y) = \lambda(Y, X)$. From the property of the Chalmers-Metcalf operator, we know that $Y$ is an invariant subspace of $T$, and thus we can calculate the trace of the restriction to $Y$ (the trace being defined in the standard way).
\end{itemize}

All of the mentioned properties and a broader discussion of the notion of the Chalmers-Metcalf functional (or operator) can be found in \cite{lewickiskrzypekoperator2}. For some applications of the Chalmers-Metcalf functional (or operator) see for example: \cite{chalmersmetcalf}, \cite{CM1}, \cite{CM2}, \cite{lewickiprophetcm}, \cite{lewickiskrzypekoperator2}, \cite{lewickiskrzypekoperator}, \cite{skrzypekmatrix}, \cite{sokolowski}. We remark, that during the process of exploring the dimension of the set minimal projection, we shall establish some new properties of this notion, which can possibly be of an independent interest (see for example Lemma \ref{lemcmdim}).

Following the definition of a minimal $I$-set, a Chalmers-Metcalf functional/operator
$$\alpha_1 x_1 \otimes f_1 + \alpha_2 x_2 \otimes f_2 + \ldots + \alpha_l x_l \otimes f_l$$
will be called minimal, if no smaller subset of $x_i \otimes f_i$ yields a Chalmers-Metcalf functional/operator for $Y$ (with any possible weights $\alpha_i$). Directly from Lemma \ref{lemosz} we get the following estimate in the real case.
\begin{lem}
\label{lemoszcm}
Let $X$ be an $n$-dimensional real normed space and let $Y \subseteq X$ be its $k$-dimensional subspace, where $1 \leq k \leq n-1$. Suppose that there exists a minimal Chalmers-Metcalf operator $T: X \to Y$ supported on $l$ points. Then $\dim \mathcal{P}_{\min}(X, Y) \leq k(n-k) - l + 1$.
\end{lem}
\emph{Proof.} Follows from Lemma \ref{lemoszcm} applied to $\mathcal{L}_Y(X, Y) \subseteq \mathcal{L}(X, X)$ and the equality $\dim \mathcal{L}_Y(X, Y)=k(n-k)$. \qed
\section{General upper bounds for the dimension of the set of minimal projections}
\label{sectgeneral}

In this section we establish upper bounds for the dimension of the set of minimal projections, which hold for every finite dimensional real normed space $X$ and its subspace $Y \subseteq X$. In the full generality, these estimates are optimal, but interestingly, it makes a certain difference if the minimal projection is of norm $1$ or not. We start with the following easy lemma.

\begin{lem}
\label{lemnorma}
Let $X$ be an $n$-dimensional real normed space. Suppose that $Y, X_0 \subseteq X$ are linear subspaces such that $Y_0 =X_0 \cap Y$ is non-empty and $\lambda(Y_0, Y)=1$. Then $\lambda(Y, X) \geq \lambda(Y_0, X_0)$
\end{lem}

\emph{Proof.} Let $Q: Y \to Y_0$ be a linear projection of norm $1$ and let $P: X \to Y$ be any linear projection. It is straightforward to check that 
$$P_0=Q \circ (P|_{X_0}):X_0 \to Y_0$$
is a linear projection. Moreover
$$\lambda(Y_0, X_0) \leq \|P_0\| = \| Q \circ (P|_{X_0}) \| \leq \|Q\| \cdot \|P\| = \|P\|,$$
which proves the inequality $\lambda(Y, X) \geq \lambda(Y_0, X_0)$. \qed

%gdzies oznaczenie normy l_1, odniesienie do operatora C-M, trace operatora

The following lemma is well-known.

\begin{lem}
\label{lemhiper}
Let $X=(\mathbb{K}^n, \| \cdot \|)$, where $\| \cdot \|= \| \cdot \|_{\infty}$ or $\| \cdot \|=\| \cdot \|_1$. Let $Y \subseteq X$ be a $(n-1)$-dimensional subspace given as $Y = \ker f$, where $f \in X^*$ is defined as $f(x) = x_1 + x_2 + \ldots + x_n$. Then $\lambda(Y, X)= 2 - \frac{2}{n}$ and the unique minimal projection $P_0:X \to Y$ is given by
$$P_0(x) = x - f(x)u,$$
where $u = \left ( \frac{1}{n}, \frac{1}{n}, \ldots, \frac{1}{n} \right ) \in \mathbb{K}^n$.
\end{lem}
\emph{Proof.} See Theorem 2.3.1 and Theorem 2.3.12 in \cite{lewickihabilitacja}. \qed

\begin{twr}
\label{twrgeneral}
Let $X$ be an $n$-dimensional real normed space and let $Y \subseteq X$ be its $k$-dimensional subspace, where $1 \leq k \leq n-1$. Then: 
\begin{enumerate}[(a)]
\item $\dim \mathcal{P}_{\min}(X, Y) \leq k(n-k)$ and this estimate is the best possible in general.
\item If $\lambda(Y, X)>1$, then $\dim \mathcal{P}_{\min}(X, Y) \leq k(n-k)-2$ and this estimate is the best possible in general. 
\end{enumerate}
\end{twr}

%\textcolor{red}{Dla hiperpłaszczyzn można wskazać przykład, że każdy wymiar pomiędzy $0$ i $n-3$ jest możliwy do osiągnięcia, ale ogólnie dla wymiaru $k$ nie wiem jak to zrobić.}

\emph{Proof.} For part (a), the upper bound is immediate, as for any fixed projection $P_0 \in \mathcal{P}(X, Y)$ we have $\mathcal{P}(X, Y) \subseteq P_0 + \mathcal{L}_Y(X, Y)$ and $\dim \mathcal{L}_Y(X, Y)=k(n-k)$. To prove that this estimate is best possible, let us consider the case $X=\ell_1^n$ and 
$$Y=\{((x_1, x_2, \ldots, x_n) \in \mathbb{R}^n \ : x_i=0 \text{ for } 1 \leq i \leq n-k\}.$$
For such $k$-dimensional subspace $Y$ we have $\lambda(Y, \ell_1^n)=1$, since for a projection $P_0: \ell_1^n \to Y$ defined as
$$P_0(x)=(0, \ldots, 0, x_{n-k+1}, \ldots, x_n).$$
it follows that, if $x \in \ell_1^n$, then
$$\|P_0(x)\| = |x_{n-k+1}| + \ldots + |x_n| \leq |x_1| + \ldots + |x_n|= \|x\|_1,$$
and in consequence $\|P_0\|=1$. For $1 \leq i \leq n-k$ we denote by $e_i^* \in Y^{\perp}$ a linear functional given by $e^*_i(x)=x_i$ and for $n-k+1 \leq i \leq n$ we denote by $e_i \in \mathbb{R}^n$, the $i$-th vector from the canonical unit basis. Let us also define 
$$L_{i, j} = e_j \otimes e_i^* \in \mathcal{L}_Y(X, Y) \quad \text{ for } 1 \leq i \leq n-k, \ n-k+1\leq j \leq n.$$
Operators $L_{i, j}$ form a basis of the space $\mathcal{L}_Y(X, Y)$, which is of dimension $k(n-k)$. To prove that $\mathcal{P}_{\min}(X, Y)=k(n-k)$, it is therefore enough to verify that for all $1 \leq i \leq n-k$, $n-k+1\leq j \leq n$ we have $P_0 + L_{i, j} \in \mathcal{P}_{\min}(X, Y)$. Thus, we need to prove that $\|P_0 + L_{i, j}\| \leq 1$. Indeed, for any $x \in \mathbb{R}^n$ we have
$$\|(P_0+L_{i, j})(x)\|=\|(0, \ldots, 0, x_{n-k+1}, \ldots, x_j + x_i, \ldots, x_n)\|$$
$$=\sum_{n-k+1 \leq m \leq n, m \neq j} |x_m| + |x_j+x_i| \leq   |x_i| + \sum_{n-k+1 \leq m \leq n} |x_m| \leq \|x\|_1,$$
which proves that $P_0 + L_{i, j}$ is also a minimal projection and the conclusion follows.

Now we shall prove the upper bound of part (b). By Lemma \ref{lemosz} it is enough to verify, that if a linear subspace $Y \subseteq X$ satisfies $\lambda(Y, X)>1$, then any Chalmers-Metcalf operator for $Y$ is supported on at least three pairs in $\ext B_X \times \ext B_{X^*}$. Let us suppose on the contrary, that there is a Chalmers-Metcalf operator $T: X \to X$ of the form
$$T= \alpha x_1 \otimes f_1 + (1-\alpha) x_2 \otimes f_2,$$
where $\alpha \in [0, 1]$ and $(x_i, f_i) \in \ext B_X \times \ext B_{X^*} $ for $i=1, 2$. In fact, we must have $\alpha \in (0, 1)$ -- if, for example, $\alpha=1$, then by taking any $y \in Y$ such that $f_1(y) \neq 0$, we would obtain that $x_1 \in Y$. However, this is clearly impossible, as $x_1$ is a norming point for a minimal projection $P: X \to Y$, that satisfies $\|P\|>1$. Let us denote $\tilde{f}_i=f_i|_Y$ for $i=1, 2$. We observe that functionals $\tilde{f}_1, \tilde{f}_2$ are linearly dependent in $Y^*$. Otherwise, we could find $y \in Y$ such that $\tilde{f}_1(y)=1$ and $\tilde{f}_2(y)=0$. In this case we would again obtain that $x_1 \in Y$, which is impossible. Thus, we can suppose that $\tilde{f}_2= c \tilde{f}_1$ for some $c \in \mathbb{R}$. Since both of these functionals are of norm $1$ in $Y^*$, we get that $c \in \{-1, 1\}$. Hence
$$T(y)=\tilde{f}_1(y)(\alpha x_1 + c(1-\alpha)x_2)$$
for $y \in Y$, which implies $\alpha x_1 + c(1-\alpha)x_2 \in Y$. Therefore, we can calculate trace of the operator $T|_Y: Y \to Y$ as
$$\tr(T|_Y) = \tilde{f}_1(\alpha x_1 + c(1-\alpha)x_2) \leq \|\tilde{f}_1\| \|\alpha x_1 + c(1-\alpha)x_2 \| \leq \alpha \|x_1\| + |c|(1-\alpha) \|x_2\| = 1.$$
Thus, from a property of the Chalmers-Metcalf operator we get
$$1<\lambda(Y, X) = \tr(T|_Y) \leq 1,$$
which is a contradiction. This proves the claimed upper bound.

To finish the proof of Theorem, we are left with showing that this upper bound is the best possible. To do this, we will use a norm, which is a certain Cartesian product of the norms $\| \cdot \|_1$ and $\| \cdot \|_{\infty}$. Clearly, there is no need to consider the case $k=1$, as $1$-dimensional subspaces are always $1$-complemented. Therefore, we assume that $2 \leq k \leq n-1$ and we define a norm $\| \cdot \|$ in $\mathbb{R}^n$ as
$$\|x\|=\max\left \{ \sum_{i=1}^{n-k+2} |x_i|, |x_{n-k+3}|, \ldots, |x_n| \right \}.$$
Thus, in the case $k=2$ we get just a $\| \cdot \|_1$ norm. We define also a $k$-dimensional subspace $Y \subseteq \mathbb{R}^n$ given by
$$Y = \{ x \in \mathbb{R}^n \ : \ x_1+x_2+x_3 =0, \ x_i=0 \text{ for } 4 \leq i \leq n-k+2\}.$$
We shall prove that $\lambda(Y, X) = \frac{4}{3}$ and $\dim \mathcal{P}_{\min}(X, Y)=k(n-k)-2.$

If $V \subseteq \ell_1^3$ is defined as 
$$V = \{ x \in \mathbb{R}^3 \ : \ x_1+x_2+x_3=0 \},$$
then by Lemma \ref{lemhiper} we have that $\lambda(V, \ell_1^3) = \frac{4}{3}$ and the unique minimal projection is given by
$$\ell_1^3 \ni x \to\left ( \frac{2x_1-x_2-x_3}{3}, \frac{-x_1+2x_2-x_3}{3}, \frac{-x_1-x_2+2x_3}{3} \right ) \in V.$$
Let $X_0 \subseteq X$ be three-dimensional subspace defined as
$$X_0=\{x \in \mathbb{R}^n \ : \ x_i = 0 \text{ for } 4 \leq i \leq n\}$$
and let $Y_0 = X_0 \cap Y$. Clearly, an operator $Q:Y \to Y_0$ defined as
$$Q(x)=(x_1, x_2, x_3, 0, \ldots, 0)$$
for $x \in Y$, is a linear projection of norm $1$. Restriction of $\| \cdot \|$ to $X_0$ yields a space linearly isometric to $\ell_1^3$. Moreover, subspaces $V$ and $Y_0$ are corresponding under this isometry. Thus, we have that $\lambda(Y_0, X_0)=\frac{4}{3}$. From Lemma \ref{lemnorma} it follows now that $\lambda(Y, X) \geq \frac{4}{3}$. It remains to observe that if we define $P_0:X \to Y$ by the formula
$$P_0(x) = \left ( \frac{2x_1-x_2-x_3}{3}, \frac{-x_1+2x_2-x_3}{3}, \frac{-x_1-x_2+2x_3}{3}, 0, \ldots, 0, x_{n-k+3}, \ldots, x_n \right ),$$
then $\|P_0\| = \frac{4}{3}$. This follows immediately from the definition of $\| \cdot \|$ and the triangle inequality
$$ \left | \frac{2x_1-x_2-x_3}{3} \right | + \left | \frac{-x_1+2x_2-x_3}{3} \right | + \left | \frac{-x_1-x_2+2x_3}{3} \right | \leq \frac{4}{3}(|x_1| + |x_2| + |x_3|) \leq \frac{4}{3} \|x\|.$$
Now we define linear functionals $g, e_4^*, \ldots, e^*_{n-k+2}$ in $Y^{\perp}$ by: $g(x)=x_1+x_2+x_3$ and $e_i^*(x) = x_i$ (where $4 \leq i \leq n-k+2$) for $x \in \mathbb{R}^n$. We also define $v, w \in Y_0$ as $v=(1, -1, 0, \ldots, 0)$ and $w=(1, 0, -1, 0, \ldots, 0)$. Finally, let $e_i \in Y$ be the $i$-th vector from the canonical unit basis for $n-k+3 \leq i \leq n$. The operators of the form $u \otimes f$, where $u \in \{v, w, e_{n-k+3}, \ldots, e_n\}$ and $f \in \{g, e_4^*, \ldots, e^*_{n-k+2}\}$ form a basis of the space $\mathcal{L}_{Y}(X, Y)$. To establish that $\dim \mathcal{P}_{\min}(X, Y) \leq k(n-k)-2$, it is enough to prove that $P_0 + \frac{1}{3} u \otimes f \in \mathcal{P}_{\min}(X, Y)$, that is
\begin{equation}
\label{norma}
\left \|P_0 + \frac{1}{3} u \otimes f \right \| \leq \frac{4}{3}
\end{equation}
for every $u \in \{v, w, e_{n-k+3}, \ldots, e_n\}$, $f \in \{g, e_4^*, \ldots, e^*_{n-k+2}\}$ and $u \otimes f \not \in \{v \otimes g, w \otimes g \}.$
In order to do this, we shall consider some cases. In the following, we assume that $x \in \mathbb{R}^n$ is a vector such that $\|x\| \leq 1$. In particular we have
\begin{equation}
\label{suma}
|x_1| + \ldots + |x_{n-k+2}| \leq 1
\end{equation}
\begin{equation}
\label{modul}
    |x_i| \leq 1 \text{ for } 1 \leq i \leq n.
\end{equation}
In each case, we shall consider only the relevant coordinates of the vector $(P_0 + \frac{1}{3} u \otimes f)(x)$. By this, we mean these coordinates, which can alter the norm of this vector, when compared with $P_0(x)$ (for which we already know that $\|P_0(x)\| \leq \frac{4}{3}$).
\begin{enumerate}
\item $f=e_{i}^*$, $u=e_j$ for $4 \leq i \leq n-k+2$ and $n-k+3 \leq j \leq n$. Then the $j$-th coordinate is equal to
$$\left (\left (P_0 + \frac{1}{3} e_j \otimes e_i^* \right )(x) \right )_j=x_j+\frac{1}{3}x_i.$$
From (\ref{modul}) and the triangle inequality we get that
$$\left | x_j+\frac{1}{3}x_i \right | \leq 1 + \frac{1}{3}= \frac{4}{3},$$
which proves (\ref{norma}).
\item $f=g$, $u=e_{j}$ for $n-k+3 \leq j \leq n$. Then the $j$-th coordinate is equal to
$$\left ((P_0 + \frac{1}{3} e_j \otimes g)(x) \right )_j=x_j+\frac{1}{3}(x_1+x_2+x_3).$$
Therefore, from the triangle inequality and conditions (\ref{suma}), (\ref{norma}) we obtain
$$\left | x_j+\frac{1}{3}(x_1+x_2+x_3) \right | \leq |x_j| + \frac{|x_1|+|x_2|+|x_3|}{3} \leq 1 + \frac{1}{3}=\frac{4}{3},$$
which proves (\ref{norma}).
\item $f=e_i^*$, $u=v$ for $4 \leq i \leq n-k+2$. In this case, we look at the sum of the absolute values of the first three coordinates.
$$\sum_{j=1}^{3} \left |\left ((P_0 + \frac{1}{3} v \otimes e_i^*)(x) \right )_j \right |$$
$$=\left | \frac{2x_1-x_2-x_3+x_i}{3} \right | + \left | \frac{-x_1+2x_2-x_3-x_i}{3} \right | + \left | \frac{-x_1-x_2+2x_3}{3} \right |$$
$$\leq \frac{4}{3}(|x_1| + |x_2| + |x_3|) + \frac{2}{3}|x_i| \leq \frac{4}{3} (|x_1| + |x_2| + |x_3| + |x_i|)$$
by the triangle inequality and (\ref{suma}). Again, this verifies inequality (\ref{norma}).
\item $f=e_i^*$, $u=w$ for $4 \leq i \leq n-k+2$. In this case, the reasoning is the same as in the previous one.
\end{enumerate}
We conclude that the set $\mathcal{P}_{\min}(X, Y)$ has the required dimension $k(n-k)-2$ and the result follows. \qed

In this paper, we focus solely on minimal projections, but we remark that the previous result could be easily generalized to the case of \emph{minimal extensions}. If $A:Y \to Y$ is a fixed linear operator, then among all possible extensions of $A$ to the linear operator from $X$ to $Y$, we can consider the extensions of minimal norm, which are called \emph{minimal extensions}. Clearly, a minimal projection corresponds to a minimal extension, when $A$ is the identity on $Y$. It is easy to verify that the set of minimal extensions of a given linear operator is convex and we can consider its dimension in the same way as for the set of minimal projections. The upper bounds from Theorem \ref{twrgeneral} hold also for minimal extensions, but the condition $\lambda(Y, X)>1$ in the second part should be replaced by the condition $\|A_0\| > \|A\|$, where $A_0$ is any minimal extension of $A$. For minimal extensions, there is also a Chalmers-Metcalf operator, which makes it possible to carry the proof of upper bound of point (b) in the exactly same manner. The upper bound of point (a) is immediate. We leave the details for the interested reader, who is referred to \cite{lewickiprophet} for the discussion on the Chalmers-Metcalf operator for minimal extensions.

%remark o minimal extensions
\section{Subspaces in a general position in polyhedral norms}
\label{sectpoly}
%%definicja general position
%%definicja norming pair
%%definicja facet
%definicja kuli i sfery

In Theorem \ref{twrgeneral} we established optimal upper bounds on the dimension of the set of minimal projections, when $Y \subseteq X$ is a linear subspace. However, it is natural to ask, to what extent these maximal values reflect what happens for a "generic" subspace $Y$ of $X$. Let us recall, that the set is called \emph{meagre}, if it is a countable union of nowhere dense sets, considered with an appropriate topology. In the case of the set of $k$-dimensional subspaces of $n$-dimensional space, one can consider the topology of Grassmannian $\Gr(k, n)$. Let us return for a moment, to the more general setting of best approximations in a linear subspace, that was considered in Section \ref{sectapprox}. In Lemma \ref{lemosz} we have established an upper bound for the dimension of the set of best approximations. However, this should be compared with the result of Zalgaller (\cite{zalgaller}), who proved that the set of $k$-dimensional subspaces $Y \subseteq X$, for which there exists $x \in X$ such that a best approximation for $x$ in $Y$ is \textbf{not} unique is meagre in the set of all $k$-dimensional subspaces, when endowed with the topology of Grassmanian. That is, it is a small subset of $\Gr(k, n)$. Speaking more precisely, Zalgaller has extented the ideas developed earlier by Ewald, Larman and Rogers (\cite{ewald}) to prove that the set of such $k$-dimensional subspaces $Y \subseteq X$ has finite $(k(n-k)-1)-$dimensional Hausdorff measure. Thus one can say that for almost all $k$-dimensional subspaces $Y$, every vector $x \in X$, has the unique best approximation in $Y$.

We should note here, that a result of Zalgaller does not seem to have any direct implications for the case of linear projections. The problem of finding a minimal projection can be considered as a problem of finding the best approximation in the subspace $\mathcal{L}_Y(X, Y)$, as explained in Section \ref{sectapprox}. However, the subspaces of the form $\mathcal{L}_Y(X, Y)$ already form a very small set in the set of all $k(n-k)$-dimensional subspaces of the space of linear operators $\mathcal{L}(X, X)$ (which is of dimension $n^2$). 

To illustrate better the case of projections, we shall consider some examples. Let us start with the case of $1$-dimensional subspaces. It is easy to see, that if $y_0 \in S_X$ is of norm $1$, then the subspace $Y=\lin\{ y_0 \}$ is $1$-complemented in $X$ and a minimal projection from $X$ onto $Y$ is unique if and only if, $y_0$ is a smooth point of $S_X$. A well-known theorem due to Mazur (\cite{mazur}) guarantees that for a finite dimensional $X$, the smooth points form a dense subset of the unit sphere of $X$. Even more, the non-smooth points form a meagre set of the unit sphere (see Theorem 5.2 in \cite{gruber}). Hence, a minimal projection is not unique only for a meagre set of $1$-dimensional subspaces.

Now, let us consider the opposite end of the spectrum. In the case of $k=n-1$, i.e. hyperplane case, there are also numerous results suggesting that the non-uniqueness of a minimal projection is a rare situation. For example, let us take $X= \ell_{\infty}^{n}$ and $Y = \ker f$, where $f=(f_1, f_2, \ldots, f_n) \in \mathbb{R}^n$ is a vector satisfying $\|f\|_1=1$. Then, if a minimal projection from $X$ onto $Y$ is not unique, we must have $f_i=0$ for some $1 \leq i \leq n$ (see Theorem 2.3.1 in \cite{lewickihabilitacja}). Hence, a minimal projection is unique on an open and dense set of hyperplanes. Similar conclusion can be derived in the case of $X = \ell_1^n$ (see Theorem 2.3.12 in \cite{lewickihabilitacja}). Soon, we shall see that this holds true generally for all polyhedral normed spaces.

For an arbitrary dimension $k$ of $Y$, the situation seems to be much more elusive. Already in the case of $k=n-2$ for $X = \ell_{\infty}^n$, the description of minimal projections is quite complicated (see \cite{lewickitwodim}, \cite{GL}, \cite{sokolowski}). The evidence seems to be pointing in the direction of non-uniqueness of minimal projection only on a meagre set, similarly like in the result of Zalgaller, but the situation is far from being clear. However, in the setting of the polyhedral norms, we are able to at least significantly improve the estimates given in Theorem \ref{twrgeneral}, for almost all subspaces. The phrase "almost all" has a certain precise meaning in our considerations. We say that two linear linear subspaces $Y, Z \subseteq \mathbb{R}^n$ are in a \emph{general position}, if $\dim(\lin{(Y \cup Z)})$ is as large as possible, that is 
$$\dim(\lin{(Y \cup Z)})=\min(\dim Y + \dim Z, n).$$

To simplify the notation further, let us suppose that $X$ is a $n$-dimensional polyhedral normed space such that $\ext B_X=\{x_1, x_2, \ldots, x_N\}$ and $\ext B_{X^*} = \{f_1, f_2, \ldots, f_M\}$. When we consider a single linear subspace $Y \subseteq X$, then we will say that it is in a \emph{general position} if it is in a general position with respect to every subspace of the form: $\lin \{x_i: i \in I \},$ where $I \subseteq \{1, 2, \ldots, N\}$ and $\bigcap_{i \in I} \ker f_i$, where $I \subseteq \{1, 2, \ldots, M\}$. It is clear, that for a fixed $1 \leq k \leq n-1$, the set of $k$-dimensional linear subspaces, which are in a general position, is an open and dense subset of all $k$-dimensional subspaces (when considered with the topology of $\Gr(k, n)$).

Troughout this section, when $X$ wil be $n$-dimensional polyhedral normed space, we will always assume that
$$\ext B_X=\{x_1, x_2, \ldots, x_N\} \: \text{ and } \: \ext B_{X^*} = \{f_1, f_2, \ldots, f_M\}.$$
With this assumption, the norming pairs for a projection $P \in \mathcal{P}(X, Y)$ are pairs of the form $(x_i, f_j)$ (where $1 \leq i \leq N$ and $1 \leq j \leq M$), for which we have $f_j(P(x_i))=\|P\|.$

We start with a lower estimate on the dimension of the set of minimal projections. The following results complements Lemma \ref{lemosz} in the setting of the polyhedral normed spaces. We remark that main idea of the proof has already appeared before. See for example Theorem 15 in \cite{sokolowski}.

\begin{lem}
\label{lemoszdol}
Let $X$ be an $n$-dimensional polyhedral normed space and let $Y \subseteq  X$ be its subspace of dimension $k$, where $1 \leq k \leq n-1$. Suppose that $P_0 \in \mathcal{P}_{\min}(X, Y)$ has at most $m$ norming pairs. Then 
$$\dim \mathcal{P}_{\min}(X, Y) \geq k(n-k)-m+1.$$
\end{lem}
\emph{Proof.} We consider a Chalmers-Metcalf functional for $Y$, which after the appropriate renaming of the vertices, is a linear functional $F$ in $\mathcal{L}^*(X, X)$ of the form
$$F(L)=\alpha_1 f_1(L(x_1)) + \ldots + \alpha_l f_l(L(x_l))$$
that satisfies $F(L)=0$ for every operator $L \in \mathcal{L}_{Y}(X, Y)$. Because a Chalmers-Metcalf functional is supported on the subset of the norming pairs of an arbitrary minimal projection, we have that $f_i(P_0(x_i))=\|P_0\|$ and $l \leq m$. For simplicity, we assume that $l=m$ and some $\alpha_i$ can be equal to $0$. Obviously, the linear functionals $(x_i \otimes f_i)|_{\mathcal{L}_{Y}(X, Y)}$ ($1 \leq i \leq m)$ are linearly dependent in the space $\mathcal{L}^*_{Y}(X, Y)$. So, let us pick a maximal linearly independent subset of this set. We can suppose that it is $(f_i \otimes x_i)|_{\mathcal{L}_{Y}(X, Y)}$ for $1 \leq i \leq d$, where $d \leq m-1$. We add 
$$N=k(n-k)-d \geq k(n-k)-m+1$$
functionals $F_1, F_2, \ldots, F_N \in \mathcal{L}^*_{Y}(X, Y)$ to complete a basis of the space $\mathcal{L}^*_{Y}(X, Y)$. For every $1 \leq j \leq N$ we can find an operator $L_{j} \in \mathcal{L}_{Y}(X, Y)$ such that:
\begin{itemize}
\item $(x_i \otimes f_i)(L_j)=0$ for every $1 \leq i \leq d$, 
\item $F_i(L_j)=0$ for every $1 \leq i \leq N$, $i \neq j$
\item $F_j(L_j)=1$. 
\end{itemize}
Because the pairs $(x_i, f_i)$, where $1 \leq i \leq m$, contain all possible norming pairs for $P_0$ we have that $(x \otimes f)(P_0) < \|P_0\|$ for every $x \in \ext B_{X}$ and $f \in \ext B_{X^*}$ such that $(x, f) \neq (x_i, f_i)$ for $1 \leq i \leq m$. Since these sets of extreme points are finite, we have that
$$C_1=\max \{(x \otimes f)(P_0) \ : \ x \in \ext B_{X}, f \in \ext B_{X^*}, (x, f) \neq (x_i, f_i) \} < \|P_0\|.$$
Let us also take $C_2 = \max \{\|L_j\| \ : 1 \leq j \leq N\}$ and define $\alpha = \frac{\|P_0\|-C_1}{C_2}$. Then, for $(x, f) \in \ext B_{X} \times \ext B_{X^*}$, $(x, f) \neq (x_i, f_i)$ and $1 \leq j \leq N$, we have that
$$(x \otimes f)(P_0+\alpha L_j) \leq C_1 + \alpha \|L_j\| \leq C_1 + \alpha C_2 = C_1 + |P_0\| - C_1 = \|P_0\|.$$
Moreover, for $1 \leq i \leq m$ we have
$$(x_i \otimes f_i)(P_0+\alpha L_j)=(x_i \otimes f_i)(P_0) + \alpha (x_i \otimes f_i)(L_j)=\|P_0\|.$$
This shows that for every $1 \leq j \leq N$ a projection $P_j \in \mathcal{P}(X, Y_0)$ defined as $P_j = P_0 + \alpha L_j$ is a minimal projection. It remains to observe, that the operators $L_1, \ldots, L_N$ are linearly independent. Indeed, if
$$a_1L_1 + a_2L_2 + \ldots + a_NL_N=0$$
for some $a_1, a_2, \ldots, a_N \in \mathbb{R}$, then
$$0=F_i(a_1L_1 + a_2L_2 \ldots + a_NL_N)=a_i,$$
for every $1 \leq i \leq N.$ This proves that $\dim \mathcal{P}_{\min}(X, Y_0) \geq N \geq k(n-k)-m+1.$ \qed

Combining the previous Lemma with Lemma \ref{lemoszcm} we immediately get the following corollary.

\begin{cor}
Let $X$ be an $n$-dimensional polyhedral normed space and let $Y \subseteq  X$ be its subspace of dimension $k$, where $1 \leq k \leq n-1$. Let $l$ be the minimal number of pairs that a Chalmers-Metcalf operator for $Y$ can have and let $m$ be the minimal number of norming pairs that a minimal projection from $X$ onto $Y$ can have. Then
$$k(n-k)-m+1 \leq \dim \mathcal{P}_{\min}(X, Y) \leq k(n-k)-l+1.$$
\end{cor}

Now, we are going to improve estimates of Theorem \ref{twrgeneral} for subspaces in a general position. As it often happens, the case of $1$-complemented subspaces requires a different approach. It turns out that for $1$-complemented subspaces in a general position the situation is quite simple, as a minimal projection is always unique.

\begin{twr}
\label{twrproj1}
Let $X$ be an $n$-dimensional polyhedral normed space. Suppose that $Y \subseteq X$ is a $k$-dimensional subspace $X$, where $1 \leq k \leq n-1$, which is in a general position and satisfies $\lambda(Y, X)=1$. Then there is a unique projection $P: X \to Y$ such that $\|P\|=1$.
\end{twr}

\emph{Proof}. Clearly, the unit ball $B_Y$ of $Y$ is a symmetric convex polytope and every face of $Y$ is contained in some face of $X$. Let $\tilde{F}$ be a facet of $B_Y$ (i.e. $(k-1)$-dimensional face) that is contained in some face $F$ of $B_X$. It is easy to see that $F$ has also to be a facet of $B_X$. Indeed, let $V \subseteq X$ be the linear span of $F$ translated to $0$ and let $\tilde{V} \subseteq Y$ be the linear span of $\tilde{F}$ translated to $0$. Then $\dim \tilde{V} = k-1$ and if $\dim V \leq n-2$, then by the assumption of general position we have
$$k-1=\dim(\tilde{V}) = \dim(Y \cap V) = \dim Y + \dim V - \dim(\lin(Y \cup V))$$
$$=k + \dim V - \min\{k+ \dim V, n\} \leq k + (n-2) - n=k-2,$$
which is a contradiction.

Because $Y^*$ is also a polyhedral space, the extreme points of its unit ball are given by the vectors corresponding to the facets of $B_Y$ (i.e. the vertices of $B_{Y^*}$). Thus, every functional on $Y$ is a convex combination of functionals corresponding to the facets of $B_Y$. In particular, we can choose functionals $\tilde{g}_1, \tilde{g}_2, \ldots, \tilde{g}_k \in \ext B_{Y^*}$, which correspond to the facets of $B_Y$ and their linear span is equal to $Y^*$. Let $g_i \in \ext B_{X^*}$ be the unique functional such that $g_i|_Y = \tilde{g}_i$ for $1 \leq i \leq k$ (where the uniqueness follows from the previous part). For $1 \leq i \leq k$ we also choose a vector $y_i \in S_Y$ satisfying
$$1=\|y_i\|=g_i(y_i)>g(y_i),$$
for all $g\in \ext B_{X^*}$, $g \neq g_i$. In other words, we choose $y_i$ from the relative interior of the facet determined by $\tilde{g}_i$. Let $P:X \to Y$ be any projection of norm $1$ and let $w \in \ker P$. For $1 \leq i \leq k$ and $\lambda \in \mathbb{R}$ with sufficiently small absolute value we have $g_i\left ( y_i + \lambda w \right ) > g \left ( y_i + \lambda w \right )$ for all $g \in \ext B_{X^*}, g \neq g_i$. Hence, $g_i$ is norming for $y_i + \lambda w$, if $\lambda$ has sufficiently small absolute value. Thus, for such $\lambda \neq 0$ we have
$$1=\|y_i\|=\left \| P\left ( y_i + \lambda w \right ) \right \| \leq \|y_i+\lambda w \| = g_i(y_i+\lambda w) = 1 + \lambda g_i(w),$$
from which we get $g_i(w)=0$ for every $1 \leq i \leq k$.

Now let us assume that there exist two projections $P_1, P_2: X \to Y$ such that $\|P_1\|=\|P_2\|=1$. We will prove that $P_1(u)=P_2(u)$ for every $u \in X$. Indeed, $u - P_1(u) \in \ker P_1$, $u-P_2(u) \in \ker P_2$ and therefore, by the previous observation we have,
$$g_i(P_1(u)-P_2(u))=g_i((u-P_1(u)))-g_i((u-P_2(u)))=0,$$
for every $1 \leq i \leq k$. Hence, because $\tilde{g}_i$ span $Y^*$, we get that $\tilde{g}(P_1(u)-P_2(u))=0$ for every $g \in Y^*$, which proves that $P_1(u)=P_2(u)$. This concludes the proof. \qed

In the case of non $1$-complemented subspaces, the following lemma is crucial for improving our estimate for the subspaces in a general position. We emphasize, that in this lemma we do not assume that $X$ is a polyhedral space or real. This result is a broad generalization of the result first proved by Cheney and Morris (\cite{cheneymorris}), who observed that if $P_0 \in \mathcal{P}_{\min}(X, Y)$ and $\|P_0\|>1$, then if from the set of norming functionals for $P_0$:
$$\{f \in \ext B_{X^*}: \ \|f \circ P_0 \| = \|P_0\|\}$$
we remove all the linear dependency relations of the form $f=-g$, then there still exists some linear dependency relation over $Y$. This was later improved by Lewicki and Skrzypek (see Theorem 2.30 in \cite{lewickiskrzypekoperator2}), who proved that the same fact holds, but when we restrict only to the functionals appearing in the Chalmers-Metcalf operator for $Y$ (which is a subset of all norming functionals). In the following lemma, we prove, that if for the functionals appearing in the Chalmers-Metcalf operator, we remove all of the linear dependency relations holding over $X$, there is still some linear dependency relation over $Y$ remaining. This can be easily stated in terms of the dimension of the subspaces spanned by the functionals.

\begin{lem}
\label{lemcmdim}
Let $X$ be an $n$-dimensional normed space and let $Y \subseteq X$ be subspace such that $\lambda(Y, X)>1$. Let $T: X \to X$ given by
$$T=\alpha_1 x_1 \otimes f_1 + \alpha_2 x_2 \otimes f_2 \ldots + \alpha_l x_l \otimes f_l,$$
be a Chalmers-Metcalf operator for $Y$. If we denote $\tilde{f}_i = f_i |_Y$ for $1 \leq i \leq l$, then the dimension of $\lin \{\tilde{f}_1, \ldots, \tilde{f}_l\}$ (considered as a subspace of $Y^*$) is strictly smaller then the dimension of $\lin \{ f_1, \ldots, f_l \}$ (considered as a subspace of $X^*$).
\end{lem}

\emph{Proof.} Clearly $\dim \{\tilde{f}_1, \ldots, \tilde{f}_l\} \leq \dim \lin \{ f_1, \ldots, f_l \}$, so let us assume that we have the equality. Suppose that
$$d = \dim \{\tilde{f}_1, \ldots, \tilde{f}_l\} = \dim \lin \{ f_1, \ldots, f_l \},$$
and $\lin \{ f_1, \ldots, f_l \} = \lin \{f_1, \ldots, f_d \}$. Then we have also that $\lin \{ \tilde{f}_1, \ldots, \tilde{f}_l \} = \lin \{\tilde{f}_1, \ldots, \tilde{f}_d \}$. For any $1 \leq i \leq l$ we have the unique representation
$$f_i = a_{i1}f_1 + \ldots + a_{id}f_d,$$
where $a_{ij} \in \mathbb{K}$. Then
$$T = \sum_{i=1}^{d} \left ( \sum_{j=1}^{l} a_{ij} \alpha_j x_j \right ) \otimes f_i.$$
By the assumption, $\tilde{f}_1, \ldots, \tilde{f}_d$ are linearly independent in $Y^*$. Hence, for every $1 \leq i \leq d$ there exists $y_i \in Y$ such that $f_i(y_i)=1$ and $f_j(y_i)=0$ for $i \neq j$. Since $T(Y) \subseteq Y$ we have
$$T(y_i) = \sum_{j=1}^{l} a_{ij} \alpha_j x_j \in Y.$$
Thus, we can calculate the trace of $T|_Y$, considered as an operator from $Y$ to $Y$, as follows
$$\tr(T|_Y) = \sum_{i=1}^{d} \tilde{f}_i \left ( \sum_{j=1}^{l} a_{ij} \alpha_j x_j \right ) = \sum_{i=1}^{d} f_i \left ( \sum_{j=1}^{l} a_{ij} \alpha_j x_j \right ) = \sum_{i=1}^{l}\alpha_i f_i(x_i).$$
We recall that by the assumption we have $\tr(T|_Y) = \lambda(Y, X)>1$. But on the other hand
$$ \left | \sum_{i=1}^{l}\alpha_i f_i(x_i) \right | \leq \sum_{i=1}^{l} \alpha_i |f_i(x_i)| \leq \sum_{i=1}^{l} \alpha_i \|f_i\| \|x_i\| = \sum_{i=1}^{l} \alpha_i=1,$$
which leads to a contradiction. The conclusion follows. \qed

As a consequence of the previous Lemma, we have the following result for non $1$-complemented subspaces in a general position in the polyhedral normed spaces.

\begin{lem}
\label{twrchalmersmetcalf}
Let $X$ be an $n$-dimensional polyhedral normed space. Let $Y \subseteq X$ be a $k$-dimensional subspace ($2 \leq k \leq n-1$), which is in a general position and satisfies $\lambda(Y, X)>1$. Then every Chalmers-Metcalf operator $T:X \to X$ for $Y$:
$$T=\alpha_1 x_1 \otimes f_1 + \alpha_2 x_2 \otimes f_2 \ldots + \alpha_l x_l \otimes f_l,$$
is supported on at least $n$ pairs, that is, $l \geq n$.

\end{lem}

\emph{Proof.} Let us assume the opposite, that there is a Chalmers-Metcal operator that is of the form
$$T=\alpha_1 x_1 \otimes f_1 + \alpha_2 x_2 \otimes f_2 \ldots + \alpha_l x_l \otimes f_l,$$
where $l \leq n-1$. Let us denote $\tilde{f}_i=f_i|_Y \in Y^*$ and let $d$ be the dimension of $\lin \{ \tilde{f}_1, \ldots, \tilde{f}_l \}$. We may further assume that $\lin \{ \tilde{f}_1, \ldots, \tilde{f}_l \}=\lin \{\tilde{f}_1, \ldots, \tilde{f}_d \}$. We will prove that $d$ is maximal possible, that is $d=k$. Indeed, let us suppose that $d \leq k-1$. We notice that if $d=l$, then the dimension of $\lin \{ \tilde{f}_1, \ldots, \tilde{f}_l \}$ in $Y^*$ would be equal to dimension of $\lin \{f_1, \ldots, f_l \}$, contradicting the previous Lemma. Therefore $d<l$. For any $d < i \leq l$ we can write
$$\tilde{f}_i = a_{i1}\tilde{f}_1 + \ldots + a_{id}\tilde{f}_d,$$
where $a_{ij} \in \mathbb{K}$. Now, let us consider a functional
$$f_i - a_{i1}f_1 - \ldots - a_{id}f_d \in X^*.$$
It vanishes on $Y$, but also on every element in the intersection $\bigcap_{j=1}^{d} \ker f_j \cap f_i$. The intersection of kernels of $d+1$ functionals in $X^*$ has the dimension at least $n-d-1 \geq n-k+1-1=n-k.$ Considering the fact that $\dim Y=k$ and the assumption of a general position of $Y$, the functional above has to vanish on all $X$ and therefore it is equal to $0$. Because this holds for every $d < i \leq l$, this proves that the dimension of $\lin \{ f_1, \ldots, f_l \}$ in $X^*$ is equal to $d$, which contradicts the previous Lemma. It follows that $d=k$.

Now, if we write
$$\tilde{f}_i = a_{i1}\tilde{f}_1 + \ldots + a_{ik}\tilde{f}_k,$$
for $k < i \leq l$, then we have
\begin{equation}
\label{kombinacja}
T=\alpha_1 x_1 \otimes \tilde{f}_1 + \ldots + \alpha_l x_l \otimes \tilde{f}_l = \sum_{i=1}^{k}  \left ( \alpha_i x_i +  \sum_{j=k+1}^{l} a_{ij} \alpha_j x_j \right ) \otimes \tilde{f}_i.
\end{equation}
Because $\tilde{f}_1, \ldots, \tilde{f}_k$ form a basis of $Y^*$, for every $1 \leq i \leq k$ we can find $y_i \in Y$ such that $\tilde{f}_i(y_i)=1$ and $\tilde{f}_j(y_i)=0$ for $i \neq j$. From the inclusion $T(Y) \subseteq Y$, we get that
$$T(y_i) = \alpha_i x_i +  \sum_{j=k+1}^{l} a_{ij} \alpha_j x_j \in Y.$$
Now we will use the assumed inequality $l \leq n-1$. Because $Y$ is in a general position with the respect to $\lin\{x_i, x_{k+1}, \ldots, x_l\}$, their intersection will be of dimension $0$ as $l-k+1 + k = l+1 \leq n$. Thus, for every $1 \leq i \leq k$ we get that
$$\alpha_i x_i +  \sum_{j=k+1}^{l} a_{ij} \alpha_j x_j=0,$$
but considering the relation (\ref{kombinacja}), this implies that $T$ vanishes on $Y$. This is clearly impossible, as $T|_Y$ has non-zero trace (equal to $\lambda(Y, X)$). This proves that $l \geq n$ and the conclusion follows. \qed

Using the previous results, we can now easily improve the estimates of Theorem \ref{twrgeneral}. As a consequence, by taking $k=n-1$, we get that a minimal projection onto a hyperplane in a general position is always unique.

\begin{twr}
\label{twrwymiar}
Let $X$ be an $n$-dimensional polyhedral normed space. Let $Y \subseteq X$ be a $k$-dimensional subspace ($2 \leq k \leq n-1$), which is in a general position. Then $$\dim \mathcal{P}_{\min}(X, Y) \leq k(n-k)-n+1.$$
\end{twr}

\emph{Proof}. By Theorem \ref{twrproj1} we can assume that $\lambda(Y, X)>1$. In this case, the conclusion follows from Lemma \ref{twrchalmersmetcalf} combined with Lemma \ref{lemoszcm}. \qed

Now we move to the question of estimating the number norming pairs of minimal projections in polyhedral norms. Because a Chalmers-Metcalf operator is supported on the norming pairs of an arbitrary minimal projection, Lemma \ref{twrchalmersmetcalf} implies that if $X$ is $n$-dimensional normed space and $Y \subseteq X$ is in general position, then \textbf{every} minimal projection from $X$ onto $Y$ has at least $n$ norming pairs. We shall establish that this holds true for every subspace (not necessarily in general position), but for \textbf{some} minimal projection. We will prove it using an approximation argument. To make it more precise, it is convenient to introduce some metric on the Grassmanian $\Gr(n, k)$, that generates its natural topology. We will use a well-known metric. For two $k$-dimensional subspaces $Y, Z$ of an $n$-dimensional real normed space $X$, we define $d(Y, Z)$ to be the Hausdorff distance of the unit spheres of $Y$ and $Z$. This metric induces the topology of the Grassmanian (see for example I-2 in \cite{ferrer}). We start with an auxiliary result, related to the continuity of the relative projection constant.

\begin{lem}
Let $X$ be finite-dimensional real normed space and let $Y \subseteq X$ be $k$-dimensional subspace (where $1 \leq k \leq n-1$). Then, there exists $C>0$ (depending on $Y$) such that for any $k$-dimensional subspace $Y_0 \subseteq X$ we have
$$\lambda(Y_0, X) \leq \lambda(Y, X) + Cd(Y, Y_0)$$
\end{lem}
\emph{Proof.} Let us fix $Y_0 \subseteq X$. We denote $d=d(Y, Y_0)$ and $\lambda=\lambda(Y, X)$. Let $y_1, y_2, \ldots, y_k  \in Y$ be unit vectors forming a basis of $Y$ and let $x_{k+1}, \ldots, x_n \in X$ be unit vectors, completing $y_i$ to a basis of $X$. By choosing $C$ large enough, we can assume that $d$ is sufficiently small, so that $\lin\left ( Y_0 \cup \{x_{k+1}, \ldots, x_n\} \right )=X$. Let us take a minimal projection $P \in \mathcal{P}_{\min}(X, Y)$. From the definition of the metric $d$, for every $k+1 \leq i \leq n$ we can find $z_i \in Y_0$ such that $\|z_i\|=\|P(x_i)\|$ and 
$$\|z_i-P(x_i)\| \leq d \|P(x_i)\| \leq d \|P\|=d \lambda.$$ 
We define a linear operator $P_0: X \to Y$ by the following conditions: 
\begin{itemize}
\item $P_0(y_0)=y_0$ for every $y_0 \in Y_0$,
\item $P_0(x_i)=z_i$ for $k+1 \leq i \leq n$.
\end{itemize}
It is clear that $P_0 \in \mathcal{P}(X, Y_0)$. We shall now estimate the norm of the operator $P_0-P$. For every unit vector $y \in Y$ we can find a unit vector $y_0 \in Y_0$ such that $\|y_0-y\| \leq d$. Hence, for a unit vector $y \in Y$ we have
$$\|P_0(y) - P(y)\| = \|P_0(y-y_0)  + (y_0-y)\| \leq \|y-y_0\|( \|P_0\| + 1) \leq d(\|P_0\| + 1).$$
Moreover, for $k+1 \leq i \leq n$ we have
$$\|P_0(x_i) - P(x_i)\| = \|z_i - P(x_i)\| \leq d \lambda.$$
Because all norms in $\mathbb{R}^n$ are equivalent, there exists $c>0$ (depending on $Y$, but not on $Y_0$) such that for all vectors $a=(a_1, a_2 \ldots, a_n) \in \mathbb{R}^n$ we have
$$\|a\|_1=\sum_{i=1}^{n} |a_i| \leq c\|a_1y_1 + \ldots + a_ky_k + a_{k+1}x_{k+1} + \ldots + a_nx_n\|.$$
Hence, for a vector $x \in \mathbb{R}^n$, given in the form above, we can combine previous estimates to obtain
$$\|(P_0-P)(x)\| = \|(P_0-P)(a_1y_1 + \ldots + a_ky_k + a_{k+1}x_{k+1} + \ldots + a_nx_n)\|$$
$$\leq \sum_{i=1}^{k} |a_i| \|(P_0-P)(y_i)\| + \sum_{i=k+1}^n |a_i| \|(P_0-P)(x_i)\|$$
$$\leq d\|a\|_1 (\|P_0\|+1 + \lambda) \leq cd\|x\| (\|P_0\| + 1 + \lambda).$$
Hence
$$\|P_0\| \leq \|P\| + \|P_0-P\| \leq \lambda + cd(\|P_0\| + 1 + \lambda),$$
and therefore
$$\|P_0\| \leq \frac{\lambda}{1-cd} + cd \frac{1+\lambda}{1-cd} \leq \lambda + Cd,$$
for some constant $C>0$ (independent of $Y_0$). \qed

\begin{lem}
\label{lemnorm}
Let $X$ be an $n$-dimensional polyhedral normed space and let $2 \leq k \leq n-1$. Let $Y \subseteq X$ be a $k$-dimensional subspace that satisfies the condition: every projection in $\mathcal{P}_{\min}(X, Y)$ has at most $m$ different norming pairs, where $m$ is a fixed positive integer. Then, there exists $r>0$ such that for every $k$-dimensional subspace $Y_0 \subseteq X$ satisfying $d(Y, Y_0) \leq r$, every projection in $\mathcal{P}_{\min}(X, Y_0)$ has at most $m$ different norming pairs.
\end{lem}

\emph{Proof}. If the thesis does not hold, we can find a sequence of $k$-dimensional subspaces $Y_i \subseteq X$ such that $d(Y, Y_i) \to 0$ and for every $i \geq 1$, there exists a minimal projection $P_i \in \mathcal{P}_{\min}(X, Y_i)$ with at least $m+1$ different norming pairs. By compactness, we can choose a subsequence of $(P_i)_{i \geq 1}$ that converges to a certain linear operator $P: X \to X$. Without loss of generality, we can suppose that $P_i \to P$ for $i \to \infty$. It is easy to see that $P$ is a projection from $X$ to $Y$. In fact, we will show that it is a minimal projection. Indeed, for every $1 \leq r \leq N, 1 \leq s \leq M$ and $i \geq 1$ we have $f_s(P_i(x_r)) \leq \lambda(Y_i, X)$, so by passing to the limit we get $f_s(P(x_r)) \leq \liminf_{i \to \infty} \lambda(Y_i, X)$. Hence
$$\lambda(Y, X) \leq \|P\| = \max_{r, s} f_s(P(x_r)) \leq \liminf_{i \to \infty} \lambda(Y_i, X).$$
Combining this with the previous Lemma, we get that $\lambda(Y_i, X) \to \lambda(Y, X)$ and it follows that $P$ is a minimal projection from $X$ onto $Y$.

Because $X$ is polyhedral, the set of possible norming pairs for any projection is finite. In particular, this set has also a finite number of $(m+1)$-element subsets and some subset will be included in the set of norming pairs for infinitely many projections $P_i$. However, if we have $f_s(P_i(x_r))=\|P_i\|$, then again by passing to the limit we obtain also that $f_s(P(x_r))=\|P\|$ (since $\|P_i\|=\lambda(Y_i, X) \to \lambda(Y, X)=\|P\|$). We conclude that norming pairs will be preserved in the limit and hence $P$ has at least $m+1$ different norming pairs. This contradicts the assumption and the proof is finished. \qed

We are ready to prove the aforementioned result, that for every subspace (not necessarily in a general position) of a polyhedral normed space, there exists a minimal projection with at least $n$ norming pairs. This may be quite surprising, as this is independent of the dimension of $Y$ and, to our best knowledge, this kind of result was not known before. It should be noted, that similar estimate does not hold in general norms. Skrzypek and Shekhtman gave an example of two-dimensional subspace of $\ell_{4}^n$, which has only six norming points and six norming functionals (see Corollary 2.6 in \cite{skrzypektwodim}).

\begin{twr}
\label{twrpunkty}
Let $X$ be an $n$-dimensional polyhedral normed space and let $Y \subseteq X$ be a $k$-dimensional subspace, where $1 \leq k \leq n-1$. Then, there exists a projection $P \in \mathcal{P}_{\min}(X, Y)$ with at least $n$ norming pairs.
\end{twr}

%\textcolor{red}{To twierdzenie pokazuje, że w normie wielościennej na każdą nie $1$-uzupełnialną podprzestrzeń istnieje projekcja z przynajmniej $n$ parami normującymi. W zasadzie dla projekcji o normie $1$ to też powinno działać, ale to trzeba pewnie trochę inaczej zapisać.}

\emph{Proof}. Let us assume that every projection in $\mathcal{P}_{\min}(X, Y)$ has at most $n-1$ norming pairs. Then by the previous Lemma, there exists $r>0$ such that for every $k$-dimensional normed space $Y_0 \subseteq X$ with $d(Y, Y_0) \leq r$ every projection in $\mathcal{P}_{\min}(X, Y)$ has at most $n-1$ norming pairs. However, we can find $Y_0$ such that $d(Y, Y_0) \leq r$ and $Y_0$ is in general position. If $\lambda(Y_0, X)=1$, then by Theorem \ref{twrproj1} we get that $\dim \mathcal{P}_{\min}(X, Y_0)=0$. On the other hand, from Lemma \ref{lemoszdol} it follows that
$$\dim \mathcal{P}_{\min}(X, Y_0) \geq k(n-k) - n+2 \geq (n-1) - n + 2 \geq 1.$$ 
This is a contradiction, which proves that $\lambda(Y_0, X)>1$. Thus, the assumptions of Lemma \ref{twrchalmersmetcalf} are satisfied and in this case, every Chalmers-Metcalf operator for $Y_0$ is supported on at least $n$ pairs. Because every pair in a Chalmers-Metcalf operator is a norming pair, this proves that every minimal projection for $Y_0$ has at least $n$ norming pairs. This is a contradiction, which concludes the proof. \qed

We do not know if the estimate given in Theorem \ref{twrwymiar} is optimal for $2 \leq k \leq n-2$, but Lemmas \ref{lemoszdol} and \ref{lemnorm} immediately imply the following: if $X$ is an $n$-dimensional polyhedral normed space and $Y \subseteq X$ is a $k$-dimensional subspace ($1 \leq k \leq n-1)$ such that every projection in $\mathcal{P}_{\min}(X, Y)$ has at most $m$ norming pairs, then the inequality 
$$\dim \mathcal{P}_{\min}(X, Y_0) \geq k(n-k)-m+1$$
holds for every $k$-dimensional subspace $Y_0$ in some open neighbourhood of $Y$. Thus, if it would be possible to have $m=n$, the estimate of Theorem \ref{twrwymiar} would turn out to be optimal. In \cite{sokolowski} (see Theorem 13 and 15) there is an example of a $k$-dimensional subspace of $\ell_{\infty}^n$ with a minimal projection that has exactly $n$ norming pairs. However, one would need to know it for every minimal projection and not only a particular one.
%\begin{twr}
%\label{twr42}
%Let $X$ be a $4$-dimensional real normed space such that $\ext B_X=\{x_1, x_2, \ldots, x_N\}$ and $\ext B_{X^*} = \{f_1, f_2, \ldots, f_M\}$. Let $Y \subseteq X$ be a $2$-dimensional subspace, which is in general position with respect to every subspace of the form $\lin \{x_i: i \in I \},$ where $I \subseteq \{1, 2, \ldots, N\}$ and $\bigcap_{i \in I} \ker f_i$, where $I \subseteq \{1, 2, \ldots, M\}$. Then there exists a unique minimal projection from $X$ onto $Y$.
%\end{twr}

%\emph{Proof}. By Theorem \ref{twrproj1} we can assume that $\lambda(Y, X)>1$.

\section{Acknowledgements}

The research of the first author was funded by the Priority Research Area SciMat under the program Excellence Initiative – Research University at the Jagiellonian University in Kraków.

\end{document}